# Identifying the global reference set in DEA: A mixed 0-1 LP formulation with an equivalent LP relaxation

*Mahmood Mehdiloozad*

*Department of Mathematics, College of Sciences, Shiraz University, Shiraz 71454, Iran*
(E-mail: m.mehdiloozad@gmail.com)

**Abstract**

The recent study by [Mehdiloozad, Mirdehghan, Sahoo, & Roshdi (2015) On the identification of the global reference set in data envelopment analysis. EJOR, 245, 779–788] proposes a linear programming (LP) model for the problem of finding the global reference set (GRS) of the evaluated decision making unit (DMU). This technical note revisits the problem and reformulates it as a mixed 0-1 LP model. By applying the LP relaxation method, it then transforms the formulated model into an equivalent LP model. Finally, it shows that the resulting LP model is equivalent to the LP model of Mehdiloozad et al. (2015).

Consider a set of $n$ observed DMUs; each uses $m$ inputs to produce $s$ outputs. Let $J$ be the index set of all the observed DMUs and let $\mathbb{R}_+^d$ be the non-negative Euclidean $d$-orthant. We denote, respectively, the input and output vectors of DMU$_j$ ($j \in J$) by $\mathbf{x}_j = (x_{1j},...,x_{mj})^T \in \mathbb{R}_+^m$ and $\mathbf{y}_j = (y_{1j},...,y_{sj})^T \in \mathbb{R}_+^s$, and the input and output matrices by $\mathbf{X} = [\mathbf{x}_1 \ ... \ \mathbf{x}_n]$ and $\mathbf{Y} = [\mathbf{y}_1 \ ... \ \mathbf{y}_n]$.

Let $o \in J$ be the index of the DMU under evaluation. Then, the RAM (range-adjusted measure) model of Cooper et al. (1999) is in the following form:

$$\rho_o = \min \quad 1 - \frac{1}{m+s}\left(\mathbf{R}^{-T}\mathbf{s}^- + \mathbf{R}^{+T}\mathbf{s}^+\right)$$

subject to
$$\mathbf{X}\lambda + \mathbf{s}^- = \mathbf{x}_o,$$
$$\mathbf{X}\lambda - \mathbf{s}^+ = \mathbf{y}_o, \qquad (1)$$
$$\mathbf{1}^T \lambda = 1,$$
$$\lambda \geq \mathbf{0},\ \mathbf{s}^- \geq \mathbf{0},\ \mathbf{s}^+ \geq \mathbf{0},$$

where the vectors $\mathbf{R}^- = (R_1^-, \ldots, R_m^-)^T$ and $\mathbf{R}^+ = (R_1^+, \ldots, R_s^+)^T$ are defined by

$$\frac{1}{R_i^-} = \max_{j \in J}\{x_{ij}\} - \min_{j \in J}\{x_{ij}\}, \quad i = 1, \ldots, m;$$
$$\frac{1}{R_r^+} = \max_{j \in J}\{y_{rj}\} - \min_{j \in J}\{y_{rj}\}, \quad r = 1, \ldots, s. \qquad (2)$$

The set of all optimal solutions of model (1) can be formulated as the set of all feasible solutions of the following system of equations:

$$\begin{bmatrix} \mathbf{X}_E & \mathbf{I}_m & \mathbf{0} \\ \mathbf{Y}_E & \mathbf{0} & -\mathbf{I}_s \\ \mathbf{1}^T & \mathbf{0} & \mathbf{0} \\ \mathbf{0} & \mathbf{R}^{-T} & \mathbf{R}^{+T} \end{bmatrix} \begin{bmatrix} \lambda \\ \mathbf{s}^- \\ \mathbf{s}^+ \end{bmatrix} = \begin{bmatrix} \mathbf{x}_o \\ \mathbf{y}_o \\ 1 \\ (m+s)(1-\rho_o) \end{bmatrix},$$
$$\begin{bmatrix} \lambda \\ \mathbf{s}^- \\ \mathbf{s}^+ \end{bmatrix} \geq \mathbf{0}, \qquad (3)$$

where $\mathbf{X}_E$ and $\mathbf{Y}_E$ are the input and output matrices of RAM-efficient DMUs, respectively.

We define $\Omega_o$ as the set of all intensity vectors that are associated with the optimal solutions of model (1). In terms of the feasible solutions of (3), $\Omega_o$ can be then expressed as

$$\Omega_o := \left\{\lambda \mid (\lambda, \mathbf{s}^-, \mathbf{s}^+) \text{ is feasible to system (3) for some slack vectors } \mathbf{s}^- \text{ and } \mathbf{s}^+\right\}. \quad (4)$$

According to Mehdiloozad et al. (2015), the global reference set (GRS) of DMU$_o$, $R_o^G$, can be found based on the following relation:

$$R_o^G = \left\{ \text{DMU}_j \mid \lambda_j^{\max} > 0 \right\}, \tag{5}$$

where $\boldsymbol{\lambda}^{\max}$ is a *maximal element* of $\Omega_o$ —an element with the maximum number of positive components. And, $\boldsymbol{\lambda}^{\max}$ can be identified by using the following LP model:

$$\begin{aligned}
&\max \quad \mathbf{1}^T \boldsymbol{\alpha} + \gamma \\
&\text{subject to} \\
&\begin{bmatrix} \mathbf{X}_E & \mathbf{I}_m & \mathbf{0} \\ \mathbf{Y}_E & \mathbf{0} & -\mathbf{I}_s \\ \mathbf{1}^T & \mathbf{0} & \mathbf{0} \\ \mathbf{0} & \mathbf{R}^{-T} & \mathbf{R}^{+T} \end{bmatrix} \begin{bmatrix} \boldsymbol{\alpha}+\boldsymbol{\beta} \\ \mathbf{s}^- \\ \mathbf{s}^+ \end{bmatrix} - \begin{bmatrix} \mathbf{x}_o \\ \mathbf{y}_o \\ 1 \\ (m+s)(1-\rho_o) \end{bmatrix} (\gamma+\nu) = \mathbf{0}, \\
&\mathbf{0} \leq \boldsymbol{\alpha} \leq \mathbf{1}, \ 0 \leq \gamma \leq 1, \\
&\mathbf{0} \leq \boldsymbol{\beta}, \ \mathbf{0} \leq \mathbf{s}^-, \ \mathbf{0} \leq \mathbf{s}^+, \ 0 \leq \nu.
\end{aligned} \tag{6}$$

Precisely, if $\left( \boldsymbol{\alpha}^*, \boldsymbol{\beta}^*, \mathbf{s}^{-*}, \mathbf{s}^{+*}, \gamma^*, \nu^* \right)$ is an optimal solution to model (6), then

$$\boldsymbol{\lambda}^{\max} = \frac{1}{1+\nu^*} \left( \boldsymbol{\alpha}^* + \boldsymbol{\beta}^* \right). \tag{7}$$

We are now going to show that model (6) can be derived from a 0-1 LP model through the LP relaxation method. In this regard, first, we develop the following mixed 0-1 LP model:[1]

$$\begin{aligned}
&\max \quad \mathbf{1}^T \boldsymbol{\alpha} + \gamma \\
&\text{subject to} \\
&\begin{bmatrix} \mathbf{X}_E & \mathbf{I}_m & \mathbf{0} \\ \mathbf{Y}_E & \mathbf{0} & -\mathbf{I}_s \\ \mathbf{1}^T & \mathbf{0} & \mathbf{0} \\ \mathbf{0} & \mathbf{R}^{-T} & \mathbf{R}^{+T} \end{bmatrix} \begin{bmatrix} \boldsymbol{\lambda} \\ \mathbf{s}^- \\ \mathbf{s}^+ \end{bmatrix} - \begin{bmatrix} \mathbf{x}_o \\ \mathbf{y}_o \\ 1 \\ (m+s)(1-\rho_o) \end{bmatrix} \delta = \mathbf{0}, \\
&\boldsymbol{\alpha} \leq \boldsymbol{\lambda}, \ \gamma \leq \delta, \\
&\boldsymbol{\alpha}: \text{binary}, \ \gamma: \text{binary}, \\
&\mathbf{0} \leq \boldsymbol{\lambda}, \ \mathbf{0} \leq \mathbf{s}^-, \ \mathbf{0} \leq \mathbf{s}^+, \ 0 \leq \delta.
\end{aligned} \tag{8}$$

To clarify the idea behind developing model (8), note that (a) $\alpha_j > 0$ ($\gamma > 0$) implies $\lambda_j > 0$ ($\delta > 0$), and (b) the vector $\boldsymbol{\alpha}$ and the variable $\gamma$ are both binary. Hence, as formally demonstrated below, model (8) identifies $\boldsymbol{\lambda}^{\max}$ by maximizing $\mathbf{1}^T \boldsymbol{\alpha} + \gamma$.

**Lemma 1** For any $\boldsymbol{\lambda} \in \Omega_o$, there exists a feasible solution $\left(\boldsymbol{\lambda}', \mathbf{s}^{-'}, \mathbf{s}^{+'}, \delta', \boldsymbol{\alpha}', \gamma'\right)$ to model (8) such that $\mathbf{1}^T \boldsymbol{\alpha}' = n^+(\boldsymbol{\lambda})$, where $n^+(\cdot)$ denotes the number of positive components of a vector.

*Proof* Let $\boldsymbol{\lambda} \in \Omega_o$. By the definition of $\Omega_o$ in (4), there exists $\mathbf{s}^-$ and $\mathbf{s}^+$ such that $\left(\boldsymbol{\lambda}, \mathbf{s}^-, \mathbf{s}^+\right)$ is a feasible solution to system (3). Then, the solution $\left(\boldsymbol{\lambda}', \mathbf{s}^{-'}, \mathbf{s}^{+'}, \delta', \boldsymbol{\alpha}', \gamma'\right)$ defined by

$$\boldsymbol{\lambda}' := \boldsymbol{\lambda}, \ \mathbf{s}^{-'} := \mathbf{s}^-, \ \mathbf{s}^{+'} := \mathbf{s}^+, \ \delta' := 1, \ \alpha'_j := \begin{cases} 1 & \lambda'_j > 0, \\ 0 & \lambda'_j = 0, \end{cases} \gamma' := 1, \tag{9}$$

is feasible to model (8) and $\mathbf{1}^T \boldsymbol{\alpha}' = n^+(\boldsymbol{\lambda})$. ∎

Based on Lemma 1, the following theorem shows that $\boldsymbol{\lambda}^{\max}$ can be found by virtue of an optimal solution of model (8).

**Theorem 1** Let $\left(\boldsymbol{\lambda}^*, \mathbf{s}^{-*}, \mathbf{s}^{+*}, \delta^*, \boldsymbol{\alpha}^*, \gamma^*\right)$ be an optimal solution to model (8). Then, $\boldsymbol{\lambda}^{\max} = \dfrac{\boldsymbol{\lambda}^*}{\delta^*}$.

*Proof* From the feasibility of system (3), it can be easily asserted that $\delta^* > 0$. Moreover, since model (8) is a maximization LP problem, $\gamma^* = 1$ and $\alpha_j^* = 1$ for any $j$ that $\lambda_j^* = 1$. This indicates that $\mathbf{1}^T \boldsymbol{\alpha}^* = n^+(\boldsymbol{\lambda}^*)$. Dividing both sides of the first set of constraints of model (10) also results that $\dfrac{\boldsymbol{\lambda}^*}{\delta^*} \in \Omega_o$, implying $\mathbf{1}^T \boldsymbol{\alpha}^* \leq n^+(\boldsymbol{\lambda}^{\max})$. Now, the equality holds immediately by Lemma 1. ∎

Theorem 1 follows that the GRS can be found with the help of model (8). However, as is known, this method is not computationally efficient when the size of the model is large. in the following, we deal effectively with this issue by demonstrating that the LP relaxation of model (8) provides an equivalent LP model, which is computationally more efficient and, hence, practically more applicable than model (8).

By relaxing the binary constraints in model (8), we transform it into the following LP model:

$$\begin{aligned}
\max \quad & \mathbf{1}^T \boldsymbol{\alpha} + \gamma \\
\text{subject to} \quad & \\
& \begin{bmatrix} \mathbf{X}_E & \mathbf{I}_m & \mathbf{0} \\ \mathbf{Y}_E & \mathbf{0} & -\mathbf{I}_s \\ \mathbf{1}^T & \mathbf{0} & \mathbf{0} \\ \mathbf{0} & \mathbf{R}^{-T} & \mathbf{R}^{+T} \end{bmatrix} \begin{bmatrix} \boldsymbol{\lambda} \\ \mathbf{s}^- \\ \mathbf{s}^+ \end{bmatrix} - \begin{bmatrix} \mathbf{x}_o \\ \mathbf{y}_o \\ 1 \\ (m+s)(1-\rho_o) \end{bmatrix} \delta = \mathbf{0}, \\
& \boldsymbol{\alpha} \leq \boldsymbol{\lambda}, \ \gamma \leq \delta, \\
& \mathbf{0} \leq \boldsymbol{\alpha} \leq \mathbf{1}, \ 0 \leq \gamma \leq 1, \\
& \mathbf{0} \leq \boldsymbol{\lambda}, \ \mathbf{0} \leq \mathbf{s}^-, \ \mathbf{0} \leq \mathbf{s}^+, \ 0 \leq \delta.
\end{aligned} \quad (10)$$

**Theorem 2** Model (10) is equivalent to model (8).

*Proof* Since model (10) is an LP relaxation of the original mixed 0-1 LP problem, the optimal objective value of model (10) is an upper bound for that of model (8). To prove the reverse, let $(\boldsymbol{\lambda}^*, \mathbf{s}^{-*}, \mathbf{s}^{+*}, \delta^*, \boldsymbol{\alpha}^*, \gamma^*)$ be an optimal solution to model (10). Then, it will suffice to show that $\gamma^* = 1$ and $\alpha_j^* = 1$ for any $j$ that $\alpha_j^* > 0$.

Since system (3) is feasible, it can be easily verified—by the way of contradiction—that $\delta^* > 0$ and, consequently, $\gamma^* > 0$. We claim that $\gamma^* = 1$. To prove our claim, assume by contradiction that $\gamma^* < 1$. By dividing both sides of the constraints of model (10) at optimality by $\gamma^*$, we have

$$\begin{aligned}
& \begin{bmatrix} \mathbf{X}_E & \mathbf{I}_m & \mathbf{0} \\ \mathbf{Y}_E & \mathbf{0} & -\mathbf{I}_s \\ \mathbf{1}^T & \mathbf{0} & \mathbf{0} \\ \mathbf{0} & \mathbf{R}^{-T} & \mathbf{R}^{+T} \end{bmatrix} \begin{bmatrix} \boldsymbol{\lambda}^*/\gamma^* \\ \mathbf{s}^{-*}/\gamma^* \\ \mathbf{s}^{+*}/\gamma^* \end{bmatrix} - \begin{bmatrix} \mathbf{x}_o \\ \mathbf{y}_o \\ 1 \\ (m+s)(1-\rho_o) \end{bmatrix} \frac{\delta^*}{\gamma^*} = \mathbf{0}, \\
& \frac{\boldsymbol{\alpha}^*}{\gamma^*} \leq \frac{\boldsymbol{\lambda}^*}{\gamma^*}, \ 1 \leq \frac{\delta^*}{\gamma^*}, \\
& \mathbf{0} \leq \boldsymbol{\alpha}^* \leq \mathbf{1}, \ 0 \leq \gamma^* \leq 1, \\
& \mathbf{0} \leq \boldsymbol{\lambda}^*, \ \mathbf{0} \leq \mathbf{s}^{-*}, \ \mathbf{0} \leq \mathbf{s}^{+*}, \ 0 \leq \delta^*.
\end{aligned} \quad (11)$$

6Based on (11), the vector $\left(\boldsymbol{\lambda}', \mathbf{s}^{-\prime}, \mathbf{s}^{+\prime}, \delta', \boldsymbol{\alpha}', \gamma'\right)$, defined by

$$\boldsymbol{\lambda}' := \boldsymbol{\lambda}^*/\gamma^*, \ \mathbf{s}^{-\prime} := \mathbf{s}^{-*}/\gamma^*, \ \mathbf{s}^{+\prime} := \mathbf{s}^{+*}/\gamma^*, \ \delta' := \frac{\delta^*}{\gamma^*}, \ \alpha'_j := \min\left\{1, \frac{\alpha_j^*}{\gamma^*}\right\}, \ \gamma' := 1, \quad (12)$$

is then a feasible solution to model (10). Since the objective function value for this solution is greater than $\mathbf{1}^T \boldsymbol{\alpha}^* + \gamma^*$, the optimality of $\left(\boldsymbol{\lambda}^*, \mathbf{s}^{-*}, \mathbf{s}^{+*}, \delta^*, \boldsymbol{\alpha}^*, \gamma^*\right)$ is contracted and our claim is hence proved.

In a similar way, it can be proved that $\alpha_j^* = 1$ for any $j$ that $\alpha_j^* > 0$ and so the proof is complete. ∎

As per Theorems 1 and 2, $\boldsymbol{\lambda}^{\max}$ can be identified by solving the relaxed LP model (10). Hence, our discussion is ended by defining $\boldsymbol{\beta} := \boldsymbol{\lambda} - \boldsymbol{\alpha}$ and $\nu := \delta - \gamma$, which imply the equivalence of model (10) and model (6).

## Notes

[1] In order to find all possible (*weakly efficient*) reference units, a different mixed 0-1 LP model can be found in Roshdi, Van de Woestyne, and Davtalab-Olyaie (2014), which is not transformed into an equivalent LP model.

## References

bibliography...bibxCooper, W. W., Park, K. S., & Pastor, J. T. (1999). RAM: A range adjusted measure of inefficiency for use with additive models and relations to other models and measures in DEA. *Journal of Productivity Analysis*, *11*, 5–42.

Mehdiloozad, M., Mirdehghan, S. M., Sahoo, B. K., & Roshdi, I. (2015). On the identification of the global reference set in data envelopment analysis. *European Journal of Operational Research*, *245*, 779–788.

Roshdi, I., Van de Woestyne, I., & Davtalab-Olyaie, M. (2014). Determining maximal reference set in data envelopment analysis. arXiv:1407.2593 [math.OC].